\documentclass[11pts]{amsart}
\usepackage{amssymb}
\usepackage{amscd}
\usepackage{graphics}
\usepackage{mathabx}
\usepackage{pdfpages}
\usepackage{mathtools}

\usepackage[
hyperref,
nonotes
]{degt}

\addtheorem{Notation}{notation}

\def\Dg:{\endgraf{\bf Dg:\enspace}\ignorespaces}
\def\Le:{\endgraf{\bf Le:\enspace}\ignorespaces}
\def\Fl:{\endgraf{\bf Fl:\enspace}\ignorespaces}

\def\sphere{\Bbb S}

\newcommand{\Zz}{\mathbb{Z}}

\newcommand{\Ss}{\mathbb{S}}

\def\CV{\Cal V}
\def\CA{\Cal A}
\def\CU{\Cal U}
\def\CCA{\CA^\circ}

\def\lk{\operatorname{\ell\mathit{k}}}
\def\vlk{\operatorname{\,\overline{\!\lk}}}
\def\sg{\QOPNAME{sg}}

\def\units{^{\times\!}}
\def\Cs{\C\units}

\def\CV{\Cal V}
\def\CA{\Cal A}
\def\CCA{\CA^\circ}

\def\ZQ<#1,#2>{\{#1,#2\}}

\def\Log{\operatorname{Log}}

\let\slope\kappa

\DeclarePairedDelimiter\floor{\lfloor}{\rfloor}

\def\vect{}

\let\defect\delta

\def\sg{\operatorname{sg}}

\def\sone{S^1\sminus1}
\let\E=E
\def\cc{_!}
\def\Milnor{\bar\mu}

\let\HH=H

\def\link#1{\href{http://katlas.math.toronto.edu/wiki/#1}{#1}}

\def\3{\color{blue}}
\def\4{\color{magenta}}
\def\6{\color{cyan}}
\def\7{\color{red}}

\title{Slopes of links and signature formulas}

\author{Alex Degtyarev}

\address{%
Department of Mathematics\\
Bilkent University\\
06800 Ankara, Turkey}

\email{degt@fen.bilkent.edu.tr}

\author{Vincent Florens}

\address{%
Laboratoire de Math\'{e}matiques et leurs applications, UMR CNRS 5142\\
Universit\'{e} de Pau et des Pays de l'Adour\\
Avenue de l'Universit\'{e}\\
BP 1155 64013 Pau Cedex, France
}

\email{vincent.florens@univ-pau.fr}

\author{Ana G.\ Lecuona}

\address{%
School of Mathematics and Statistics\\
University of Glasgow\\
G12 8QQ Glasgow, Scotland, UK}

\email{ana.lecuona@glasgow.ac.uk}

\thanks{%
The first author was partially supported by the T\"{U}B\DOTaccent{I}TAK grant 118F413.
The second author was partially supported by the ANR Project LISA}



\begin{document}

\begin{abstract}
We present a new invariant, called \emph{slope}, of a colored link in an
integral homology sphere and use this invariant to complete the signature
formula for the splice of two links. We develop a number of ways of computing
the slope and a few vanishing results.
Besides, we discuss the concordance invariance of the
slope and establish its close relation to the Conway polynomials, on the one
hand, and to the Kojima--Yamasaki $\eta$-function (in the univariate case)
and Cochran invariants, on the other hand.
\end{abstract}

\maketitle

\section{Introduction}\label{S.intro}

 The principal goal of this short note
is to present a brief account of our recent results and announce a few further
advances concerning the new invariant of colored links, called
\emph{slope}.
Originally introduced as a means of handling an extra correction term in our
signature formula (\cf. \autoref{th.signature} \vs. \autoref{th.signature2}),
the slope proved to be interesting on its own right. It is closely related
to, but not strictly dependent on, the Conway potentials of the links involved
(\cf. \autoref{th.Conway} \vs. \autoref{ex.Conway}) and, as such, it can be
regarded as a generalization, to both multicomponent links and those with
nonvanishing linking numbers, of the Kojima $\eta$-function (see \autoref{coro.eta}).
In the case of two-component links,
we observe a similar relation to Milnor's $\Milnor$-invariants (see
\autoref{prop.mu=slope} \vs. \autoref{ex.equal.mu}); we expect
to generalize these results to more components. The slope is easily
computable in many ways, see \autoref{s.Fox}, \autoref{th.Conway}, or
\autoref{S.C-complex} (new results).
The most important new aspect is the concordance
invariance of the slope, see \autoref{S.concordance}.

The principal object of our study is
a \emph{$\mu$-colored link},
\ie, an oriented link~$L$ in an integral homology sphere~$\Ss$ equipped with a
surjective map $\pi_0(L)\onto\{1,\ldots,\mu\}$, called
\emph{coloring}.
The union of the components of~$L$ given the same color $i=1,\ldots,\mu$ is
denoted by $L_i$.
We denote by $X:=\Ss\sminus T_L$
the complement of a small open
tubular neighborhood of $L$.
The group $H_1(X)$ is free abelian, generated by the classes $m_C$
of the meridians of the components $C \subset L$. The coloring induces an epimorphism
\[
\varphi: \pi_1(X) \onto H:=\bigoplus_{i=1}^\mu \Z t_i
\label{eq.H}
\]
sending $m_C$ to $t_i$ whenever $C \subset L_i$.
A multiplicative character $\Go\:\pi_1(X)\to\C\units$
is determined by its values on the meridians, and the torus of
characters preserving the coloring
(\ie, those that factor through~$\Gf$)
is naturally identified with $(\C\units)^\mu$.
Often,
we split the components of~$L$ into two groups,
$L=L'\cup L''$, on which the coloring takes, respectively, $\mu'$ and $\mu''$
values, $\mu'+\mu''=\mu$; then, we regard
a character as a ``vector''
$\omega=(\omega',\omega'')\in(\C\units)^{\mu'}\times(\C\units)^{\mu''}$,
\cf. \autoref{s.signature.formula}.

Given a topological space~$X$ and a multiplicative character
$\Go\:\pi_1(X)\to\C\units$, we denote by $H_*(X;\C(\Go))$ the homology
of~$X$ with
coefficient in the local system $\C(\Go)$ twisted by~$\Go$.

\section{The signature formula}\label{S.signature}

\subsection{Signature of colored links}

Let $L$ be a $\mu$-colored link in an integral homology sphere~$\Ss$.
A \emph{spanning pair} for $(\Ss,L)$ is a pair $(N,F)$, where
$N$ is a compact smooth oriented $4$-manifold such that $\partial N=\Ss$ and
$F= F_1 \cup \ldots \cup F_\mu\subset N$ is a properly immersed surface (see~\cite{DFL2})
such that
$\partial F_i = F_i \cap \partial N=L_i$ for all $i=1,\dots,\mu$.
Fix an open tubular neighborhood $T_F$
of~$F$ and let $W_F:=N\sminus  T_F$.
We assume that the group $H_1(N\sminus F)$
is freely generated by the
meridians of the components of~$F$.
Then,
any character~$\Go$
on~$\Ss\sminus L$
extends to a unique character on $W_F$, also denoted by~$\Go$.
If $\Go\in(S^1)^\mu$ is unitary, the space $H_2(W_F;\C(\Go))$ has a Hermitian
intersection form, and we denote by $\sign^\Go(W_F)$ its signature.

\definition\label{def.signature}
The \emph{signature} and \emph{nullity}
of a $\mu$-colored link $L\subset\Ss$ are the
maps
\[*
\aligned
\Gs_L\:(S^1\sminus1)^\mu&\longrightarrow\Z,\\
\Go&\longmapsto\sign^\Go(W_F) - \sign(W_F),
\endaligned
\qquad
\aligned
n_L\:(S^1\sminus1)^\mu&\longrightarrow\Z,\\
\Go&\longmapsto\dim H_1(\Ss\sminus L;\C(\Go)),
\endaligned
\]
respectively.
By convention, we extend these functions to
all unitary characters $\Go\in(S^1)^\mu$
by patching the components~$L_i$ on which
$\Go_i=1$.
According to~\cite{DFL2} or~\cite{Viro:twisted}, $\Gs_L$ is independent of the choice
of a spanning pair.
\enddefinition


\subsection{Splice of colored links} \label{Ss.splice}
A \emph{$(1,\mu)$-colored link} is a $(1+\mu)$-colored link of the form
\[*
K\cup L=K\cup L_1\cup\ldots\cup L_\mu,
\]
where the \emph{knot}~$K$ is the only component given the distinguished
color~$0$.
We define the \emph{linking
vector} $\vlk(K,L)=(\lambda_1,\ldots,\lambda_\mu)\in\Zz^\mu$, where
$\lambda_i:=\lk(K,L_i)$ and $\lk$ is the linking number in the integral homology sphere~$\Ss$.
In the following definition, for a $(1,\mu^*)$-colored link
$K^*\cup L^*\subset\Ss^*$, $*=\prime$ or $\prime\prime$, we denote by
$T_{K^*}\subset\Ss^*$ a
small tubular neighborhood of~$K^*$ disjoint from~$L^*$. Let
$m^*,\ell^*\subset\partial T_{K^*}$ be, respectively,
its meridian and longitude in the homology sphere $\Ss^*$.

\definition \label{def:splice}
Given two $(1,\mu^*)$-colored links $K^*\cup L^*\subset\Ss^*$,
$*=\prime$ or $\prime\prime$, their \emph{splice}
is the $(\mu'+\mu'')$-colored link $L'\cup L''$ in the integral
homology sphere
$$
\Ss:=(\Ss'\sminus\operatorname{int}T_{K'})\cup_\varphi(\Ss''\sminus\operatorname{int}T_{K''}),
$$
where the gluing homeomorphism $\varphi\colon\partial T_{K'}\to\partial T_{K''}$
takes~$m'$ and $\ell'$ to~$\ell''$ and~$m''$, respectively.
\enddefinition

\subsection{(Non)-additivity of the signature}\label{s.signature.formula}

The \emph{index} of $x\in\R$ is
$\ind(x):=\lfloor x\rfloor-\lfloor-x\rfloor\in\ZZ$.
The \emph{$\Log$-function}
$\Log\colon S^1\to[0,1)$ sends $\exp(2\pi it)$
to $t\in[0,1)$.
For an integral vector $\vect\lambda\in\Zz^\mu$, $\mu\ge0$, we define the \emph{defect function}
\begin{align*}
\textstyle
\defect_\lambda \colon(S^1)^\mu&\longrightarrow\Z \\
 \vect\omega&\longmapsto\textstyle
 \ind\bigl(\sum_{i=1}^\mu \lambda_i \Log\omega_i\bigr)-\sum_{i=1}^\mu \lambda_i \ind(\Log\omega_i).
\end{align*}

\theorem[see \cite{DFL}]\label{th.signature}
For $*=\prime$ or $\prime\prime$, consider a $(1,\mu^*)$-colored link
$K^*\cup L^*\subset\Ss^*$, and let $L\subset\Ss$ be the splice of the two
links.
Let $\lambda^*:=\vlk(K^*,L^*)$ and,
for characters $\omega^*\in(S^1)^{\mu^*}$, denote
$ \upsilon^*:=(\omega^*)^{\lambda^*}\in S^1$.
Then, \emph{assuming that $(\upsilon',\upsilon'')\ne(1,1)$}, one
has
\[*
\begin{aligned}
\sigma_L(\omega',\omega'')&=
 \sigma_{K'\cup L'}(\upsilon'',\omega')+
 \sigma_{K''\cup L''}(\upsilon',\omega'')
 +\defect_{\lambda'}(\omega')\defect_{\lambda''}(\omega''),\\
n_L(\omega',\omega'')&=
 n_{K'\cup L'}(\upsilon'',\omega')+
 n_{K''\cup L''}(\upsilon',\omega'').
\end{aligned}
\]
\endtheorem

 In the special case
$\upsilon'=\upsilon''=1$, the formulas of \autoref{th.signature} are
no longer valid.
If $\upsilon^*=1$, the character $\Go^*$ is \emph{admissible}
(see \autoref{def.admissible} below)
and, hence, there is a well-defined \emph{slope}
$\slope^*:=(K^*\!/L^*)(\Go^*)$ in $\R\cup\infty$ (see \autoref{go}).
These two slopes give rise to extra correction terms shown in \autoref{fig.main}.
Here,
the domain of $\Delta\sigma(\slope',\slope'')$ and
$\Delta n(\slope',\slope'')$ is the square
$[-\infty,\infty]^2\ni(\slope',\slope'')$, and
the curve in the figures is the hyperbola $\slope'\slope''=1$.
For $\Delta\Gs$, there also is an ``explicit'' formula
\[*
\Delta\Gs(\slope',\slope'')=\sg\slope'-\sg\left(\frac{1}{\slope'}-\slope''\right),
\]
where
\[*
\sg x=\begin{cases}
 \hphantom{-}0,&\text{if $x=0$ or $\infty$},\\
 \hphantom{-}1,&\text{if $x>0$},\\
 -1,&\text{if $x<0$},
\end{cases}
\]
and we disambiguate $\infty-\infty=0$.
It is worth mentioning that, in the next statement, the terms
$\Delta\sigma$ and $\Delta n$ are the \emph{only} contribution of the
knots $K^*$; the rest
depends on the links $L^*$ only.

\theorem[see \cite{DFL2}]\label{th.signature2}
For $*=\prime$ or $\prime\prime$, consider a $(1,\mu^*)$-colored link
$K^*\cup L^*\subset\Ss^*$, and let $L\subset\Ss$ be the splice of the two
links.
Let \smash{$\Go^*\in(S^1)^{\mu^*}\!$} be two admissible
characters
\rom(so that $\upsilon'=\upsilon''=1$\rom),
and denote
$\slope^*=(K^*\!/L^*)(\omega^*)$. Then,
with
$\Delta\Gs,\Delta n\in\{0,\pm1,\pm2\}$
as in \autoref{fig.main},
\figure
\includegraphics[width=9cm]{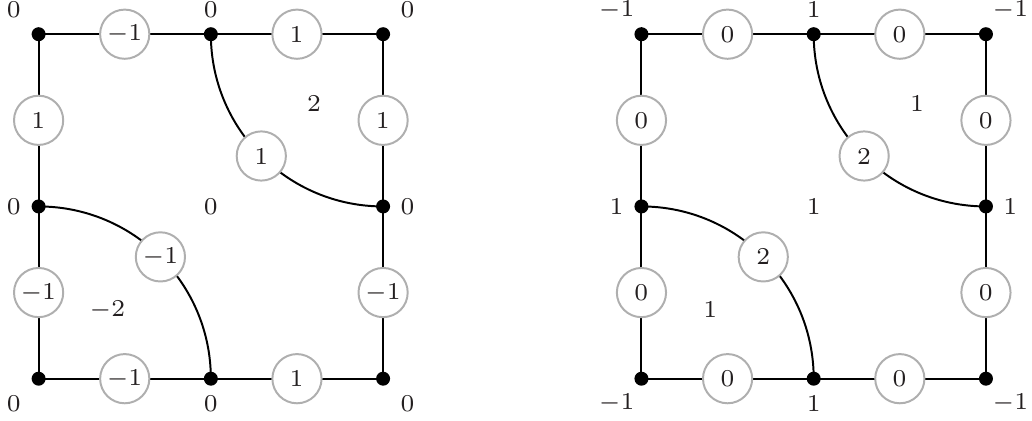}
\caption{The correction terms
$\Delta\sigma$ (left) and $\Delta n$ (right)
in \autoref{th.main}.}
\label{fig.main}
\endfigure
one has
\[*
\begin{aligned}
\sigma_{L}(\Go',\Go'')&
 =\sigma_{L'}(\omega')+\sigma_{L''}(\omega'')
 +\delta_{\lambda'}(\Go') \delta_{\lambda''}(\Go'')
 +\Delta\sigma(\slope',\slope''),\\
n_{L}(\Go',\Go'')&
 =n_{L'}(\omega')+n_{L''}(\omega'')
 +\Delta n(\slope',\slope'').
\end{aligned}
\]
\endtheorem

\section{The slope}\label{S.slope}
In this section we define the slope of a $(1,\mu)$-colored link: it is a
function defined on a subset of the character torus and taking values in the
Riemann sphere $\Cp1=\C\cup\infty$.

\subsection{Definition of the slope}\label{s.def.slope}
Fix a $(1,\mu)$-colored link $K\cup L\subset\Ss$ and let
$X:=\Ss\sminus T_L$ and $\bar X = \Ss\sminus T_{K \cup L}$.
We denote by~$m$ and~$\ell$, respectively, the meridian of $K$ and its
\emph{Seifert longitute}, \viz. the one unlinked with~$K$.
We orient~$m$ so that $m\circ\ell=1$
\emph{with respect to the orientation of~$\partial T_K$ induced from~$T_K$}.


\definition \label{def.admissible}
A character $\Go : \pi_1( X)\longrightarrow \C^\times$ is
called \emph{admissible} if $\Go([K])=1$.
\enddefinition
The variety of admissible characters is denoted
$$
\CA(K/L) = \bigl\{ \Go \in (\Cs)^\mu\bigm|\Go^\Gl=1\bigr\},
$$
where \smash{$\Gl:=\vlk(K,L)$} is the linking vector.
If $\Gl=0$, we have $\CA(K/L)=(\Cs)^\mu$;
otherwise,
letting $N:=\gcd(\Gl)$ and
$\nu:=\Gl/N$, the irreducible over~$\Q$ components of $\CA(K/L)$ are
the zero sets $\CA_d(K/L)$
of the cyclotomic polynomials $\Phi_d(\Go^\nu)$, $d\divides|N$.
We will mainly deal with \emph{nonvanishing}  admissible
characters
$$\Go\in\CCA(K/L):=\CA(K/L)\cap(\Cs\sminus1)^\mu.$$
Any character $\Go\in\CCA(K/L)$ induces a character
$\Go\:\pi_1(\bar X)\to\Cs$ sending $m$ and~$\ell$ to~$1$.




Consider the intersection $\partial_K \bar X=\partial T_K$ of
$\partial \bar X$ with the closure of $T_K$ and  the inclusion
$$ i : \partial_K \bar X \hookrightarrow \bar X.$$
Denote $Z(\Go):=\ker i_*$, where
$$ i_* : H_1(\partial_K \bar X; \C(\Go)) \longrightarrow H_1(\bar X;\C(\Go)).$$
If $\Go\in\CCA(K/L)$ is admissible,
$m$ and $\ell$ form a basis of $H_1(\partial_K \bar X;\C(\Go))=\C^2$.

\definition\label{def.slope}
Let $\Go \in \CCA(K/L)$ and assume that $\dim Z(\Go)=1$, so that $Z(\Go)$ is
generated by a single vector $am+b\ell$ for some $[a:b] \in \Cp1$.
The \emph{slope} of $K\cup L$ at~$\Go\in\CCA(K/L)$ is the quotient
\[*
(K/L)(\Go):=-\frac{a}{b}\in\Cp1=\C\cup\infty.
\]
As in \autoref{def.signature}, we extend the slope to
the whole set $\CA(K/L)$ by patching
the components of~$L$ on which $\Go$ vanishes.
\enddefinition

Thus, strictly speaking, the slope is defined on a \emph{subset} of
$\CA(K/L)$. This subset is dense (see \autoref{th.properties} below);
moreover, according to the next statement, the slope is always well defined
and \emph{real} on the \emph{unitary} admissible characters, so that
\autoref{th.signature2} makes sense.

\proposition[see~\cite{DFL2}]\label{go}
Let $\Go \in \CU(K/L):=\CA(K/L)\cap(S^1 \sminus 1)^\mu $ be a
\emph{unitary} character.
Then $\dim Z(\Go)=1$
and, hence, $(K/L)(\Go)$ is well defined. Moreover,
$(K/L)(\Go) \in \mathbb{R} \cup \infty$.
\endproposition

\proposition[essentially, \cite{DFL2}]\label{prop.properties}
For a $(1,\mu)$-colored link $K\cup L\subset\Ss$
one has\rom:
\roster
\item\label{i.orientation}
$({-K}/L)(\Go)=(K/{-L})(\Go)=(K/L)(\Go)=(K/L)(\Go\1)$\rom;
\item\label{i.component}
$(K/L')(\Go_1,\Go_2,\ldots)=(K/L)(\Go_1\1,\Go_2,\ldots)$,
where $L':={-L_1}\cup L_2\cup\ldots\cup L_\mu$\rom;
\item\label{i.mirror}
$(\bar K/\bar L)=-(K/L)$, where $\bar K\cup\bar L$ is the mirror of $K\cup L$\rom;
thus, if $K\cup L$ is amphichiral, then
$K/L$ takes values in $\{0,\infty\}$\rom;
\item\label{i.ball}
if $K$ lies in a ball disjoint from~$L$, then $K/L\equiv0$.

\item\label{i.Seifert}
if $K$ bounds a Seifert surface disjoint from a $C$-complex for~$L$,
then $K/L\equiv0$\rom;
\endroster
\endproposition

\subsection{Fox calculus}\label{s.Fox}
The slope can be computed by means
of Fox calculus from a presentation of the fundamental group
$\pi_1(\bar X)$ of the link complement,
together with the classes $m,\ell\in\pi_1(\bar X)$.
In the case of links in~$S^3$, both pieces of data can be
derived from the link diagram.
Indeed, for the group one can choose the Wirtinger
presentation,
where
meridians are the generators.
For~$\ell$, we trace a curve $C$ parallel to $K$ and such that $\lk(K,C)=0$;
then, starting from the segment corresponding to the chosen meridian of~$K$ and
moving along~$C$ in the positive direction,
we write down the corresponding generator (or its inverse)
each time when undercrossing positively
(respectively, negatively) the diagram of $K\cup L$.
Thus, let
\[*
m,\ell\in\pi_1(\bar X)=\bigl<x_1,\ldots,x_p\bigm|r_1,\ldots,r_q\bigr>
\]
and, using the epimorphisms
$F:=\<x_1,\ldots,x_p\>\onto\pi_1(\bar X)\onto\pi_1(X)\onto H$,
\cf.~\eqref{eq.H},
specialize the Fox derivatives to
$\partial/\partial x_i\:F\to\Lambda:=\Z\HH$.
Consider the complex of
$\Lambda$-modules \[* S_*\:\quad S_2\overset{\partial_1}\longto
S_1\overset{\partial_0}\longto
 S_0\longto 0,
\]
where
\[*
S_2=\bigoplus_{i=1}^q\Lambda r_i,\quad
S_1=\bigoplus_{i=1}^p\Lambda dx_i,\quad
S_0=\Lambda
\]
and $dx_i$ stands for a formal generator corresponding to~$x_i$.
The ``differential'' of a word $w\in F$ is
\[*
dw:=\sum_{i=1}^p\frac{\partial w}{\partial x_i}dx_i\in S_1;
\]
then, letting
\[*
\partial_1\:r_i\mapsto dr_i,\qquad
\partial_0\:dx_i\mapsto(\text{the image of $x_i$ in $\HH\subset\Lambda$})-1,
\]
we obtain a complex computing the homology $H_{\le1}$ of the $\HH$-covering
of~$\bar X$.

Now, pick an admissible nonvanishing character $\omega \in  \CCA(K/L)$ and
consider the specialization $S_*(\Go):=S_*\otimes_\Lambda\C(\Go)$.
Then, it is straightforward that
\[*
Z(\Go)=\Ker\bigl[ H_1(\partial_K \bar X;\CC(\omega))= \CC m \oplus \CC \ell
 \overset{\inj_*}\longto
 S_1(\Go) / \Im \partial_1(\Go)\bigr],
\]
where the inclusion homomorphism $\inj_*$ is the specialization of
$m\mapsto dm$, $\ell\mapsto d\ell$. (Note that, by the assumption that
$\omega\in\CCA(K/L)$, this homomorphism lands into $\Ker\partial_0(\Go)$.)
Computing the above kernel in the basis $m,\ell$, we can also compute the slope
whenever it is defined.

\example[the Whitehead link, see \cite{DFL2}]\label{ex.whitehead}
Consider the $(1,1)$-colored Whitehead link $K \cup L$,
see \link{L5a1} in \cite{KAT}. Since $\lk(K,L)=0$, we have
$\CCA(K/L)=\sone$.
The ``standard'' presentation of
$\pi_1(\bar X)$ (derived from the Wirtinger presentation) is
\[*
\pi_1(\bar X)=\bigl\<m, m_1,\ell \bigm|[m,\ell]=1, \ell=m_1 m\1m_1\1m m_1\1m\1m_1m \bigr\>,
\]
where $m$ and $m_1$ are meridians of $K$ and $L$, respectively,
and $\ell$ is a Seifert longitude of~$K$.
We can further specialize $\Lambda$ to the group ring
$\ZZ H_1(\bar X)=\Z[t^{\pm1},t_1^{\pm1}]$, sending
$m\mapsto t$, $m_1\mapsto t_1$, $\ell\mapsto1$.
Then, denoting by $x$, $y$ the two relations in the presentation above, we have
\[*
dx= (t-1)d\ell, \quad
dy= d\ell - t^{-1}(1-t_1)(1-t_1^{-1})dm - (1-t^{-1})(1-t_1^{-1})dm_1.
\]
The specialization at a character $\omega \in  \CCA(K/L)$ means sending
$t_1\mapsto\omega$ and  $t\mapsto1$, so that the image $\Im\partial_1(\Go)$
is generated by $dy\mapsto d\ell-(1-\omega)(1-\omega^{-1})dm$.
Thus,
$$ (K/L)(\omega)= (1-\omega)(1-\omega^{-1}).$$
In particular,  notice that the slope is not invariant under link homotopies: the Whitehead link is link homotopic to the unlink of two components, which has trivial slope.

\endexample

\subsection{The rationality}
The
\emph{characteristic varieties} $\CV_r(X)$ of~$X$ (related to $\Gf$)
are defined \via
\[*
\CV_r(X):=\bigl\{\Go\in(\Cs)^\mu\bigm|\dim H_1(X;\C(\Go))\ge r\bigr\},
 \quad r\ge0.
\]
They are algebraic varieties in $(\Cs)^\mu$,
which
are nested ($\CV_r\supset\CV_{r+1}$) and
depend  only on the fundamental group $\pi$ of~$X$ (and $\Gf$).
The irreducible components of $\CV_r(X)$ of codimension $\le1$ are closely related to the $\Z \HH$-module~$H_1(X;\Z \HH)$.
 They constitute the
zero locus of the
\emph{$(r-1)$-st order}
\[*
\Delta_{X,r-1}:=\gcd E_{r-1}(H_1(X;\Z\HH))
 =\gcd E_{r}(S_1/\Im\partial_1)\in\Z \HH,
\]
where  $E_s(M)\subset\Z \HH$ is the $s$-th elementary ideal.
 In particular, the $0$-th order $\Delta_X:=\Delta_{X,0}$ is called the
\emph{multivariate Alexander polynomial} of~$X$ and $\Gf$.
As usual, if $X$ and $\Gf$ are as in~\eqref{eq.H}, we abbreviate
$\Delta_{L,r}:=\Delta_{X,r}$ and $\Delta_L:=\Delta_X$.

\theorem[see \cite{DFL2}]\label{th.properties}
Pick a component $\CA\subset\CA(K/L)$
and let $r$ be the minimal integer such
that $\Delta_{L,r}|_{\CA}\ne0$, \ie, $\CA\sminus\CV_{r+1}(L)$ is dense
in~$\CA$.
Denote by $R$
the coordinate ring of~$\CA$
and fix a normalization of~$\Delta_{L,r}$.
Then, either
\roster
\item\label{rational.finite}
there exists a unique
polynomial $\Delta_{\CA}\in R$ such that
\[*
(K/L)(\Go)=\frac{\Delta_{\CA}(\Go)}{\Delta_{L,r}(\Go)}
\]
holds for each character $\Go\in\CCA\sminus\CV_{r+1}(L)$,
or
\item\label{rational.infinite}
the slope $(K/L)(\Go)=\infty$ is well defined
and infinite
at each character~$\Go$ in a
certain dense Zariski open subset of $\CA$.
\endroster
Case \rom{\iref{rational.infinite}} cannot occur if $r=0$, \ie, if
$\Delta_L|_\CA\ne0$, \cf. \autoref{th.Conway} below.
\endtheorem

According to this theorem, the slope gives rise to a rational function,
possibly identical $\infty$, on the variety of admissible characters, \ie, an
element of $\Q(\CA(K/L))\cup\infty$. Note that this function does not depend
on the coloring too much: one can start with the maximal coloring
(assigning
its own color to each component), upon which it suffices to identify the
variables (components of $\Go$) corresponding to the components of~$L$ given
the same color.

In the next theorem, we evaluate the Conway
potential at the radical \smash{$\sqrt\Go:=(\sqrt{\Go_1},\ldots,\sqrt{\Go_\mu})$},
which is not quite well defined. We use the convention that one of the values
of each radical is chosen and used consistently \emph{throughout the whole formula}.
The nature of the formula guarantees that the result is independent of
the initial choice.

\theorem[see \cite{DFL2}]\label{th.Conway}
For a $(1,\mu)$-colored link $K\cup L\subset\sphere$, denote
\smash{$\displaystyle\nabla':=\frac\partial{\partial t}\nabla_{K\cup L}$},
where $t$ is the first variable, corresponding to~$K$.
Then, for a character $\Go\in\CA(K/L)$, one has
\[*
(K/L)(\Go)=-\frac{\nabla'(1,\sqrt\Go)}{2\nabla_L(\sqrt\Go)}\in\C\cup\infty,
\]
provided that the expression in the right hand side makes sense, \ie,
$\nabla'(1,\sqrt\Go)$ and $\nabla_L(\sqrt\Go)$ do not vanish simultaneously.
In particular, the slope is well defined in this case.
\endtheorem

The $\eta$-function was defined by Kojima and Yamasaki \cite{KY} for
two-component links in $S^3$ with linking number zero, generalizing Goldsmith's
invariants \cite{Goldsmith}.
The next corollary follows from Jin~\cite{Jin} computing the $\eta$-function in terms of the Alexander polynomials of $K \cup L$ and $L$.

\corollary[see \cite{DFL2}] \label{coro.eta}
For any two-component link $K\cup L\subset S^3$,
$\lk(K,L)=0$, and
character $\Go \in\CA(K/L)$ such that $\Delta_{L}(\omega)\neq 0$, one has
 $$ (K/L)(\Go) = \eta_{K \cup L} (\Go).$$
\endcorollary

As an example, for the $(1,1)$-colored Whitehead link we have
$\nabla_{K \cup L}(t,t_1)= (t-t^{-1})(t_1 - t_1^{-1})$ and
$\nabla_L(t_{1})=1/(t_{1}-t_{1}^{-1})$.
Hence, for any
$\Go\in\C\units$,
$$
(K/L)({\omega})= -\bigl(\sqrt{\omega_1}-\sqrt{\omega_1}^{-1}\bigr)^2=(1-\omega_1)(1-\omega_1^{-1}),
$$
which agrees with \autoref{ex.whitehead}.
This example illustrates also the independence of the ratio
in \autoref{th.Conway} of the choice of~$\sqrt\Go$.
\autoref{th.Conway} is inconclusive if
$\nabla_L(\sqrt\Go)=\nabla'(1,\sqrt\Go)=0$, \cf. the next example.
Furthermore, even in the univariate case and at an isolated common root,
l'H\^{o}pital's rule does \emph{not} apply, \cf. \autoref{ex.Hopital} below.

\example[equal higher orders, see \cite{DFL2}]\label{ex.Conway}
Let $K \cup L_1 \cup L_2$ and $K' \cup L'_1 \cup L'_2$ be the links \link{L11n353}
and \link{L11n384} (see~\cite{KAT}),
respectively.
Both have $11$ crossings and $3$ components, and all their orders are
equal:
$$
\Delta_{K \cup L}= \Delta_{K' \cup L'}= (t_2-1)(t-1)^3(t_1-1), \quad
\Delta_{L}= \Delta_{L'}= 0,\quad
\Delta_{L,1}=\Delta_{L',1}=1.
$$
Note that \autoref{th.Conway} is inconclusive.
Since $\vlk(K,L)=\vlk(K',L')=(0,0)$, one has
$\CCA(K/L)=\CCA(K'/L')=(\sone)^2$,
and, by Fox calculus (from the link diagrams),
for any character $\Go:=(\Go_1,\Go_2)\in(\C\units)^2$,
we obtain
\[*
(K/L)(\omega)=-\frac{\omega_1 \omega_2^2 + \omega_1^2 - 4 \omega_1 \omega_2 + \omega_2^2+\omega_1}{\Go_1\Go_2},
\quad
(K'/L')(\omega)=-\frac{(\omega_1 -1)(\omega_1 \omega_2^2-1)}{\Go_1\Go_2}.
\]
Thus,
the slope can distinguish links with equal Alexander polynomials and higher
orders.

\endexample

\section{Slopes via $C$-complexes}\label{S.C-complex}

A $C$-complex (see~\cite{Cimasoni.Florens}) for a $\mu$-colored link
$L=L_1\cup\ldots\cup L_\mu\subset S^3$ is a collection $F=\bigcup_iF_i$ of Seifert
surfaces~$F_i$, $\partial F_i=L_i$, disjoint except for a finite number of
\emph{clasps}, \ie, disjoint segments connecting pairs of points $a_i\in L_i$,
$a_j\in L_j$, $i\ne j$ (and disjoint from~$L$ otherwise) along which $F_i$
and~$F_j$ intersect transversally. Each homology class in $H_1(F)$ can be
represented by a collection of \emph{proper loops}, \ie,
lops $\Ga\:S^1\to F$ such that the pull-back of each clasp is a single
segment (possibly empty). We routinely identify classes, loops, and their
images.

Given a vector $\Ge\in\{\pm1\}^\mu$, the \emph{push-off} $\Ga^\Ge$ of a
proper loop~$\Ga$ is the loop in $S^3\sminus F$ obtained by a slight shift
of~$\Ga$ off each surface~$F_i$ in the direction of~$\Ge_i$ (if $\Ga$ runs
along a clasp $F_i\cap F_j$, the shift respects both
directions~$\Ge_i$ and~$\Ge_j$). Due to~\cite{Cimasoni.Florens}
(notice that our conventions differ), this
operation gives rise to a well-defined homomorphism
$\Theta^\Ge\:H_1(F)\to H_1(S^3\sminus F)$ and we can define the
\emph{Seifert forms}
\[*
\theta^\Ge\:H_1(F)\otimes H_1(F)\to\Z,\qquad
\Ga\otimes\Gb\mapsto\lk(\Ga,\Gb^\Ge).
\]
Now, given a character $\Go\in(\C\units\sminus1)^\mu$, we define
\[*
\Pi(\Go):=\prod_{i=1}^\mu(1-\Go_i)\in\C,\qquad
A(\Go):=\sum_{\Ge\in\{\pm1\}^\mu}
 \prod_{i=1}^\mu\Ge_i\Go_i^{(1-\Ge_i)/2}\Theta^\Ge
 \:H_1(F;\C)\to H^1(F;\C).
\]
Let also
\[*
\E(\Go):=\Pi(\Go)\1A(\Go)\:H_1(F;\C)\to H^1(F;\C);
\]
if $\Go\in(S^1\sminus1)^\mu$ is unitary, this differs from $H(\Go)$
considered in~\cite{Cimasoni.Florens} by
the positive real factor $\Pi(\bar\Go)\1\Pi(\Go)\1$.
This operator computes~$\nabla_L$ (see~\cite{CiTo}), $\Gs_L$, $n_L$ and the Blanchfield pairing (see~\cite{ConFriTof}).

\theorem[see~\cite{Cimasoni.Florens}]\label{th.Cimasoni.Florens}
If
$\Go\in(S^1\sminus1)^\mu$ is a nonvanishing unitary character, then the
operator $\E(\Go)$
is Hermitian and one has
$\Gs_L(\Go)=\sign\E(\Go)$ and
$n_L(\Go)=\dim\Ker\E(\Go)+b_0(F)-1$.
\endtheorem

Now, let $K\cup L$ be a $(1,\mu)$-colored link, $\vlk(K,L)=0$. Then,
we have a well-defined class
\[*
\kappa\in H^1(F ;\mathbb{C}),\qquad \kappa\:\Ga\mapsto\lk(\Ga, K).
\]

\theorem[to appear]\label{th.C-complex}
In the notation introduced above, for $\Go\in(\Cs\sminus1)^\mu$, one has
\[*
(K/L)(\Go)=\begin{cases}
 -\<\Ga,\kappa\>, & \mbox{if } \kappa\in\Im\E(\Go)\cap\Ker\E(\Go)^\perp,\\
 \infty, & \mbox{if } \kappa\notin\Im\E(\Go)\cup\Ker\E(\Go)^\perp,\\
 \mbox{\rm undefined}, & \mbox{otherwise},
\end{cases}
\]
where, in the first case, $\Ga\in H_1(F)$ is any class such that
$\E(\Go)(\Ga)=\kappa$.
\endtheorem

\corollary\label{cor.Seifert}
Let $K\cup L\subset S^3$ be a $(1,1)$-colored link, $\lk(K,L)=0$. Pick a
Seifert surface~$F$ for~$L$ disjoint from~$K$ and consider the
Seifert
form~$\theta:=\theta^{+}$, the associated map
$\Theta\:H_1(F)\to H^1(F)$, and the linking coefficient
$\kappa\in H^1(F;\mathbb{C})$ with~$K$. For $\Go\in\C\units\sminus1$, let
$A(\Go):=\Theta-\Go\Theta^*$. Then
\[*
(K/L)(\Go)=\begin{cases}
 -(1-\Go)\<\Ga,\kappa\>, & \mbox{if } \kappa\in\Im A(\Go)\cap\Ker A(\Go)^\perp,\\
 \infty, & \mbox{if } \kappa\notin\Im A(\Go)\cup\Ker A(\Go)^\perp,\\
 \mbox{\rm undefined}, & \mbox{otherwise},
\end{cases}
\]
where, in the first case, $\Ga\in H_1(F;\C)$ is any class such that
$A(\Go)(\Ga)=\kappa$.
\endcorollary

The homomorphism $\Theta\:H_1(F)\to H^1(F)$ associated to
the Seifert form has the property that
\[
\mbox{all nonzero invariant factors of $\Theta-\Theta^*$ are $\pm1$}.
\label{eq.theta}
\]
(Indeed, $\Theta-\Theta^*$ is the intersection index, which is unimodular
modulo boundary.) Conversely, any matrix satisfying~\eqref{eq.theta} can be
realized as the Seifert form of a certain link~$L$ \cite{levine70}. Since
also any
homomorphism $H_1(F)\to\Z$ can obviously be realized by a knot algebraically unlinked
with~$L$, we have the following corollary, describing all univariate slopes
in purely algebraic terms.

\corollary\label{cor.algebraic}
Let $G$ be a finitely generated free abelian group, $\Theta\:G\to G^\vee$ a
homomorphism satisfying~\eqref{eq.theta}, and $\kappa\in G^\vee$.
Denote $A(\Go):=\Theta-\Go\Theta^*\:G\otimes\C\to G\otimes\C$.
Then, there exists a $(1,1)$-colored link whose
slope at $\Go\in\C\units\sminus1$ is given by
\[*
(K/L)(\Go)=\begin{cases}
 -(1-\Go)\<\Ga,\kappa\>, & \mbox{if } \kappa\in\Im A(\Go)\cap\Ker A(\Go)^\perp,\\
 \infty, & \mbox{if } \kappa\notin\Im A(\Go)\cup\Ker A(\Go)^\perp,\\
 \mbox{\rm undefined}, & \mbox{otherwise},
\end{cases}
\]
where, in the first case, $\Ga\in G\otimes\C$ is any class such that
$A(\Go)(\Ga)=\kappa$.
\endcorollary

Choosing a standard geometric basis
for the unimodular (modulo kernel)
skew-symmetric bilinear form $\theta\:G\otimes G\to\Z$ associated to~$\Theta$ in
\autoref{cor.algebraic}, one can restate~\eqref{eq.theta} as follows: there
is a constant $r\le\frac12\rank G$ such that the matrix~$\theta$ satisfies
the condition
\[*
\theta_{ij}-\theta_{ji}=\begin{cases}
 1, & \mbox{if } j=i-1\in\{1,3,\ldots,2r-1\}, \\
 0, & \mbox{for all other pairs $j\le i$}.
\end{cases}
\]
In this case, one can find a link~$L$ with $(\rank G-2r+1)$ components.
Using this observation and \autoref{cor.algebraic}, one can easily construct
examples of univariate slopes with prescribed properties.

\example[l'H\^{o}pital's rule, \cf. \autoref{th.Conway}]\label{ex.Hopital}
Fix integers $b\ne0,-1$, $c\ne0$, $x$, $y$ and let
\[*
\theta= \bmatrix0&b\\b+1&c\endbmatrix,\qquad
\kappa=[x,y]^t.
\]
By \autoref{cor.algebraic}, there exists a two-component link $K\cup L$ whose
slope is given by
\[*
(K/L)(\Go)=-\frac{(\Go-1)^2x(2by+y-cx)}{(\Go b-b-1)(\Go b+\Go-b)}
\]
for all $\Go\ne\Go_\pm:=(1+1/b)^{\pm1}$. Letting
$\kappa_+:=[2b+1,c]^t$ and $\kappa_-:=[0,1]^t$, we have:
\roster*
\item
if $\kappa=\kappa_\pm$, then $K/L\equiv0$ except that $(K/L)(\Go_\pm)$ is
undefined;
\item
if $\kappa=0$, then $K/L\equiv0$;
\item
for all other values of~$\kappa$, we have
$(K/L)(\Go_\pm)=\infty=\lim_{\Go\to\Go_\pm}(K/L)(\Go)$.
\endroster
Now, take for~$\theta$ the direct sum of two copies of the above matrix and
let $\kappa:=\kappa_+\oplus\kappa_-$. For the new link $K\cup L$ we have
$K/L\equiv0$ away from~$\Go_\pm$, whereas
$(K/L)(\Go_\pm)=\infty\ne\lim_{\Go\to\Go_\pm}(K/L)(\Go)$,
\ie, l'H\^{o}pital's rule does not apply in \autoref{th.Conway}.
\endexample

\section{Concordance invariance}\label{S.concordance}


This section contains most of our new results. Proofs will appear elsewhere.

\subsection{Slopes and concordance}

Two $(1,\mu)$-colored links $K \cup L$ and $K' \cup L'$
 in the same integral homology sphere~$\Ss$
are said to be \emph{concordant}
if they bound a union $D \cup A=D \cup A_1 \cup \dots \cup A_\mu$
of pairwise disjoint locally flat cylinders embedded in
$\Ss\times I$,
where $I=[0,1]$.
More precisely,
$D$ is a single cylinder with $\partial D=(K \times 0) \sqcup ({-K} \times 1)$
and
each $A_i$, $i=1,\dots,\mu$, is a disjoint union of cylinders such that
$\partial A_i=(L_i \times 0) \sqcup ({-L_i} \times 1)$.

Let
\[*
U_\mu := \bigl\{ p \in  \Z[t_1^{\pm 1},\dots, t_\mu^{\pm 1}]  \bigm| p(1,\dots,1)= \pm 1 \bigr\}.
\]
A character
$\Go\in(\Cs)^\mu$ is called a \emph{concordance root} if
$p(\Go)=0$ for some $p\in U_\mu$ (\emph{cf}.\ the \emph{Knotennullstelle}~\cite{NagelPowell}).
We denote by $\CA\cc(K/L)  \subset \CA(K/L)$ the
subset of the characters which are \emph{not} concordance roots.

For example, for $\mu=1$, a root of unity $\Go$ is \emph{not} a concordance root
if and only if the order of~$\Go$ is a prime power (including $\Go=1$ of
order~$1$). It follows that the \emph{closure} $\bar\CA\cc(K/L)$ is equal to
$\CA(K/L)$ if $\vlk(K,L)=0$; otherwise, $\bar\CA\cc(K/L)$ is the union of the
$\Q$-components $\CA_d(K/L)$ corresponding to the divisors
$d\divides|\gcd\vlk(K,L)$ that are prime powers (including $d=1$), see~\autoref{s.def.slope}.

\theorem \label{th.main}
Let $K \cup L$ and $K' \cup L'$ be two concordant $(1,\mu)$-colored links.
Then, there is a natural identification $\CA\cc(K/L)=\CA\cc(K'/L')$ and,
for all $\Go \in \CA\cc(K/L)$,
$$ (K/L)(\Go)= (K'/L')(\Go).$$
\endtheorem

\corollary\label{cor.concordance}
If two $(1,\mu)$-colored links $K\cup L$ and $K'\cup L'$ are concordant,
then
the slopes $K/L$ and $K'/L'$,
\emph{restricted to $\Q(\bar\CA\cc(K/L))\cup\infty$}
\rom(\cf. \autoref{th.properties}\rom), coincide.
\endcorollary

\example\label{ex.KL.vs.LK}
Note that, for a
maximally colored
$\mu$-component link, \autoref{cor.concordance} gives us
$\mu$ \latin{a priori} distinct concordance invariants, as one can take for~$K$
any of the $\mu$ components.

For example, for the two-component link \link{L10n2}
in~\cite{KAT} we have that $K/L\ne L/K$. Indeed, using \autoref{th.Conway} we compute $K/L\equiv0$, whereas
\[*
(L/K)(\Go)=-\frac{(\Go-1)^4}{\Go^4-3\Go^3+5\Go^2-3\Go+1}.
\]
It follows, in particular, that this link is not slice, \ie,
it is not concordant to the unlink of two components.

\endexample

\subsection{Two-component links in the sphere}
Let
$K\cup L\subset S^3$ be a two-component link with
$\lk(K,L)=0$.
The
Sato--Levine invariant of this link is defined from the intersection of
Seifert surfaces of $K$ and $L$~\cite{Sato}. Cochran extended this
construction and defined a sequence of concordance invariants $\beta^i$ as
the Sato--Levine invariants of successive
derivations of $K \cup L$~\cite{Co}.
He proved that the $\beta^i$ are
certain well-defined canonical lifts of Milnor's linking numbers of the form
$\Milnor(1^{2i}2^2)$. Moreover, he showed that the sequence of $\beta^i$ is equivalent to the $\eta$-function up to a certain
change of variable~\cite[Theorem 7.1]{Co};
thus, Corollaries \ref{coro.eta} and \ref{cor.Seifert}
show that the slope
$K/L\in\Q(\C\units)$ can be used to compute
Cochran's $\Gb^i$.
Note that the relation between $\Gb^i$ and the slope can be shown directly
from the computation with Seifert surfaces, see \autoref{cor.Seifert}.


As far as the $\Milnor$ Milnor invariants are concerned, note also that they \latin{per se} do \emph{not} determine $K/L$
(\cf. \autoref{ex.equal.mu} below);
however, the \emph{vanishing} of all $\Milnor$-invariants does imply the
vanishing of the slope away from the concordance roots
(in fact, from the roots of $\Delta_L\in U_1$).
Indeed, in this case the
$\Milnor$ are well defined in~$\Z$
and the above ``lift'' to $\Gb^i$ is redundant.
\autoref{ex.KL.vs.LK} shows that
the converse does not hold: $K/L\equiv0$ does not imply
$\Milnor\equiv0$.

\proposition[see {\cite[Theorem 5.4]{Co90}} and \autoref{th.Conway}]\label{prop.mu=slope}
Let
$K \cup L\subset S^3$ be a two-component link such that all
$\Milnor(1^m2^n)$
vanish.
Then $(K/L)(\omega)=0$ for
all $\omega \in \Cs\sminus\Delta_L\1(0)\supset\CA\cc(K/L)$.
\endproposition

\example[equal $\Milnor$-invariants, \cf.~\cite{KY}]\label{ex.equal.mu}
Let $K \cup L_n\subset S^3$ be the family of links
in \autoref{fig.equal.mu}, $n \geq 2$.
\figure
\cpic{KLn}
\caption{The links $K \cup L_n\subset S^3$ in \autoref{ex.equal.mu}.}\label{fig.equal.mu}
\endfigure
They all have the same
$\Milnor$-invariants, see \cite[Application 3]{KY}.
Using the ``obvious'' Seifert surface, we immediately obtain
$\kappa=[1,0]$ and
\[*
\theta=\bmatrix k&0\\1&-1\endbmatrix\quad\mbox{for $n=2k$ even},\qquad
\theta=\bmatrix k+1&0\\1&1\endbmatrix\quad\mbox{for $n=2k+1$ odd}.
\]
Hence, by \autoref{cor.Seifert}, the slopes
\[*
(K/L_n)(\Go)=-\frac{(\Go-1)^2}
 {\floor*{\frac{n+1}{2}}\Go^2 - \bigl(2\floor*{\frac{n}{2}}+1\bigr)\Go+ \floor*{\frac{n+1}{2}}}
\]
are pairwise distinct, which agrees with the computation of~$\eta$
in~\cite{KY}, \cf. \autoref{coro.eta}.

\endexample

In conclusion, it is worth emphasizing
that, unlike
the $\eta$-function, the slope can be defined and distinguish links
even if the two components $K$ and~$L$ are  algebraically linked, that is, if $\lk(K,L)\ne0$.

\example[linked components]\label{ex.linked}
Let $K\cup L$ and $K'\cup L'$ be, respectively, the links \link{L4a1} and
\link{L7n1} in~\cite{KAT}.
Both have $\lk(K,L)=2$ (so that $\CCA=\CCA\cc=\{-1\}$) and
equal $\Milnor$-invariants:
\[*
\Milnor(12)=2,\quad \Milnor(112)=\Milnor(122)=1 \bmod 2,\quad \text{and}\quad
\Milnor(i_1 i_2 \dots i_n)=0\quad \mbox{for $n\ge4$},
\]
see~\cite{Mu}.
A simple computation using Fox calculus (see \autoref{s.Fox}) shows that
\[*
(K/L)(-1)=(L/K)(-1)=-2\quad\mbox{whereas}\quad
(K'/L')(-1)=\frac{2}{3},\quad
(L'/K')(-1)=6.
\]
Since $-1$ is not a concordance root, the two links are not concordant.
\endexample

\let\.\DOTaccent
\def\cprime{$'$}
\bibliographystyle{amsplain}
\bibliography{bibli}

\def\cprime{$'$}
\providecommand{\bysame}{\leavevmode\hbox to3em{\hrulefill}\thinspace}
\providecommand{\MR}{\relax\ifhmode\unskip\space\fi MR }
\providecommand{\MRhref}[2]{%
  \href{http://www.ams.org/mathscinet-getitem?mr=#1}{#2}
}
\providecommand{\href}[2]{#2}
\begin{thebibliography}{10}

\bibitem{KAT}
Dror Bar-Natan, Scott Morrison, and et~al., \emph{The {K}not {A}tlas}, \href
  {http://katlas.org} {\texttt{http://katlas.org}}.

\bibitem{CiTo}
D.~Cimasoni, \emph{A geometric construction of the {C}onway potential
  function}, Comment. Math. Helv. \textbf{79} (2004), no.~1, 124--146.
  \MR{2031702}

\bibitem{Cimasoni.Florens}
David Cimasoni and Vincent Florens, \emph{Generalized {S}eifert surfaces and
  signatures of colored links}, Trans. Amer. Math. Soc. \textbf{360} (2008),
  no.~3, 1223--1264 (electronic). \MR{2357695 (2009b:57009)}

\bibitem{Co}
Tim~D. Cochran, \emph{Geometric invariants of link cobordism}, Comment. Math.
  Helv. \textbf{60} (1985), no.~2, 291--311. \MR{800009}

\bibitem{Co90}
\bysame, \emph{Derivatives of links: {M}ilnor's concordance invariants and
  {M}assey's products}, Mem. Amer. Math. Soc. \textbf{84} (1990), no.~427,
  x+73. \MR{1042041}

\bibitem{ConFriTof}
Anthony Conway, Stefan Friedl, and Enrico Toffoli, \emph{The {B}lanchfield
  pairing of colored links}, Indiana Univ. Math. J. \textbf{67} (2018), no.~6,
  2151--2180. \MR{3900365}

\bibitem{DFL}
Alex Degtyarev, Vincent Florens, and Ana~G. Lecuona, \emph{The signature of a
  splice}, Int. Math. Res. Not. IMRN (2017), no.~8, 2249--2283. \MR{3658197}

\bibitem{DFL2}
\bysame, \emph{Slopes and signatures of links}, To appear, \arXiv{1802.01836},
  2018.

\bibitem{Goldsmith}
D.~L. Goldsmith, \emph{A linking invariant of classical link concordance}, Knot
  theory ({P}roc. {S}em., {P}lans-sur-{B}ex, 1977), Lecture Notes in Math.,
  vol. 685, Springer, Berlin, 1978, pp.~135--170. \MR{521732}

\bibitem{Jin}
G.~T. Jin, \emph{On {K}ojima's {$\eta$}-function of links}, Differential
  topology ({S}iegen, 1987), Lecture Notes in Math., vol. 1350, Springer,
  Berlin, 1988, pp.~14--30. \MR{979331}

\bibitem{KY}
S.~Kojima and M.~Yamasaki, \emph{Some new invariants of links}, Invent. Math.
  \textbf{54} (1979), no.~3, 213--228. \MR{553219}

\bibitem{levine70}
J.~Levine, \emph{An algebraic classification of some knots of codimension two},
  Comment. Math. Helv. \textbf{45} (1970), 185--198. \MR{266226}

\bibitem{Mu}
Kunio Murasugi, \emph{On {M}ilnor's invariant for links. {II}. {T}he {C}hen
  group}, Trans. Amer. Math. Soc. \textbf{148} (1970), 41--61. \MR{259890}

\bibitem{NagelPowell}
Matthias Nagel and Mark Powell, \emph{Concordance invariance of
  {L}evine-{T}ristram signatures of links}, Doc. Math. \textbf{22} (2017),
  25--43. \MR{3609203}

\bibitem{Sato}
N.~Sato, \emph{Cobordisms of semiboundary links}, Topology Appl. \textbf{18}
  (1984), no.~2-3, 225--234. \MR{769293}

\bibitem{Viro:twisted}
O.~Viro, \emph{Twisted acyclicity of a circle and signatures of a link}, J.
  Knot Theory Ramifications \textbf{18} (2009), no.~6, 729--755. \MR{2542693}

\end{thebibliography}

\end{document}